\documentclass[11pt,a4paper]{amsart}
\usepackage[latin1]{inputenc}
\usepackage[mathscr]{eucal} 
\usepackage{amsmath}
\usepackage{amsfonts}
\usepackage{amssymb} 
\newtheorem{theorem}{Theorem}
\newtheorem{example}[theorem]{Example}
\newtheorem{proposition}{Proposition}
\newtheorem{definition}{Definition}

\def\down #1{ I_{#1}}
\def\up #1{ F_{#1}}
\begin{document}

\title[Poset valuations and metrics]{Valuations and  metrics on partially ordered sets}
\author[C. Orum, C. Joslyn]{Chris  Orum, Cliff A. Joslyn}
\address{Chris Orum \\ Department of Mathematics, Oregon State University, Corvallis, OR 97331}
\email{orum@math.utah.edu}
\address{Cliff Joslyn \\ Pacific Northwest National Laboratory \\PO Box 999 MS K7-28, Richland, WA 99352 USA}
\email{cliff.joslyn@pnl.gov}


\maketitle
\begin{abstract}
We extend the definitions of   upper and lower valuations on partially ordered sets, and consider  the metrics they induce, in particular the metrics available (or not) based on the logarithms of such valuations.   Motivating applications in computational linguistics and computational biology are indicated.  
\end{abstract}




\section{Introduction}
%
This expository  note is motivated by our answer, given herein as Propositions \ref{prop:logs.val} and \ref{prop:logs.val.new}, to the following question:  let $P = (P, \leq)$ be a poset with  an upper or lower valuation $v(x) : P \rightarrow \bf{R}^+$; then is $\ell(x) = \log v(x)$ necessarily an upper or lower  valuation? (These terms are defined below.)  

The question arises from the common practice in information systems
(see e.g.\ \cite{BuAHiG06}) of using measures of ``semantic
similarity'' in large taxonomic vocabularies such as
WordNet\footnote{http://wordnet.princeton.edu} (in computational
linguistics) or 
the Gene Ontology\footnote{http://geneontology.org} (in computational biology) \cite{LoPStR03b}.  
Such similarity measures  
are based on a quantification of
information content as ${\mathcal I}(x) = -\log p(x )$, where $p(x) $ is a kind of 
cumulative probability defined on a poset $P$
representing the hierarchical structure of the taxonomy.  As such, $p(x)$
has the form stated in Proposition \ref{prop:a} and is often a lower
valuation.

This  question  also  arises   from Example \ref{ex:group.theory}, that deals with valuations on the $\wedge$-semilattice $\mathscr L$ of finite subgroups $X$ of a given group $G$: both $c(X) = |X|$ and $v(X) = \log c(X)$ are lower valuations on $\mathscr L$.  
Notwithstanding this example,   the logarithm of a (positive) lower valuation need not be  a lower valuation.   On the other hand, the logarithm of a positive upper valuation is always an upper valuation.

By focusing on this question   we bring together some results  (some of which are only implicit in \cite{MoB81}, the primary predecessor of this work)   
concerning the practical differences between upper and
lower valuations defined on partially ordered sets; and we describe the 
metrics they induce.    We extend the previous definitions of upper and lower valuations to allow for antitone maps (instead of requiring that valuation be isotone).  The symmetry introduced by this extension allows us to consider the composition $\log ( K \cdot v(x) + A ) $ where  $v(x)$ is an upper valuation or a lower valuation, $K \in {\mathbf R}/\{0\}$,  $A \in \mathbf{R}$, and $K \cdot v(x) + A > 0$. 
  
Distance formulas involving $\mathcal I (x) = -\log p(x)$  appear in the literature.     We note  that such a formula introduced by Jiang and Conrath \cite{JiJCoD97} does not, in general, define a metric on a partially ordered set.  (Under the tacit assumption that the poset is a tree, however,  it does yield a metric.)

 While the literature on lattice valuations extends back to  Wilcox and Smiley (1939) \cite{Wilcox:Smiley:39,Wilcox:Smiley:46}  and  Birkhoff (1940) \cite{Birkhoff:79},
the literature on general poset valuations is quite thin: we are only aware of  
\cite{MoB81}, \cite{Ramana.Murty:Engelbert:85}.

\section{Preliminaries and notation} 
In the sequel $P = (P, \leq)$ always denotes a partially ordered set.   
The  \emph{greatest lower bound} or \emph{meet}  of two elements $x, y \in P$  
need not exist, but if it does it is denoted $x \wedge y$.
An ordered set $P$ in which  $x \wedge y$   always exists is  $\wedge$-semilattice.        
The \emph{least upper bound} or \emph{join}  of two elements $x, y \in P$  need not exist, but if it does it is denoted $x \vee y$.        An ordered set $P$ in which  $x \vee y$   always exists is  $\vee$-semilattice.  If $P$ is both a $\wedge$-semilattice and $\vee$-semilattice then $P$ is a lattice.  
   
We write $a \prec c$ if $c$ covers $a$ ($a \leq b \leq c$ and $a \neq c$  implies $a = b$ or $b = c$). The notation $\{a, b\} \prec \{ c, d\}$ means both $c$ and $d$ cover both $a$ and $b$, etc.   Given a subset $S \subseteq P$,   
$\mathbf{min}(S) \subseteq S$ and $\mathbf{max}(S) \subseteq S $  denote the minimal and maximal  elements of $S$ respectively.  
If $P$ has a  \emph{unique} minimal or maximal element, it is  denoted by 0 or 1, respectively.    If  $0 \in P$ and $1 \in P$ then $P$ is said to be a bounded.

Given an element $x \in P$, $\up x  = \{x' \in P : x \leq x'\}$ is the \emph{principal filter} generated by $x$. 
Given an element $x \in P$, $\down x  = \{x' \in P : x' \leq x \}$ is the \emph{principal ideal} generated by $x$.  

We shall call a finite $\vee$-semilattice $L$ in which $\down x \cap \down y = \varnothing$ for all $x,y \in L$ a \emph{tree}. 

If  $S$ and $T$ are  nonempty subsets of a multiplicative group $G= (G, \cdot , e) $, then $ST = \{ st : s \in S \text{ and }  t \in T\}$.  The index of a subgroup  $S \subseteq G$ is $[G\! :\! S] $, and if $G$ is finite $[G\!:\! S]= |G| / |S|$.  If $S$ and $T$  are subgroups of $G$ the smallest subgroup containing both $S$ and $T$ is denoted $S \vee T$.
\section{Valuations and metrics} 
Let $P$ be a poset.  
A function $f:P \rightarrow \mathbf R$ is \emph{isotone} if $x \leq y$ implies $f(x) \leq f(y)$ and \emph{strictly isotone} if $x < y$ implies $f(x) < f(y)$.    
It is \emph{antitone} if $x \leq y$ implies $f(x) \geq f(y)$ and \emph{strictly antitone} if $x < y$ implies $f(x) >  f(y)$.  
Assuming $f$ is monotone  (that is, either isotone or antitone) we use the notation 
\begin{align}
f^-(x,y) &= 
  \begin{cases}
             \sup \{ f(z) : z \in \down x \cap \down y \}, 
             & \text{ if $f$ is isotone},  \\[1ex]   
             \inf \{ f(z) : z \in \down x \cap \down y \}, 
             & \text{ if $f$ is antitone}, 
  \end{cases} \\[1ex]
f^+(x,y) &= 
  \begin{cases}
            \inf \{ f(z) : z \in \up x \cap \up y \}, 
            & \text{ if $f$ is isotone},  \\[1ex]   
            \sup \{ f(z) : z \in \up x \cap \up y \}, 
            & \text{ if $f$ is antitone}. 
  \end{cases}  
\end{align}
Note that $\down x \cap \down y$ or $\up x \cap \up y$  may be empty.  We use the convention $\inf \varnothing = +\infty$  and  $\sup \varnothing = -\infty$.   
\begin{definition}\label{def:1}     
Let $P$ be a poset.   
An isotone (antitone) function $v: P \rightarrow \mathbf R$ is a \emph{lower valuation}
if for all $x,y \in P$,  $\down x \cap \down y \neq \varnothing$ ($\up x \cap \up y \neq \varnothing$) and  
\begin{equation}\label{eq:def.lower}
v(x) + v(y)  \leq v^-(x,y) + v^+(x,y).
\end{equation} 
An isotone (antitone) function $v: P \rightarrow \mathbf R$ is an  \emph{upper valuation}
if for all $x,y \in P$, $\up x \cap \up y \neq \varnothing$ ($\down  x \cap \down y \neq \varnothing$) and 
\begin{equation}\label{eq:def.upper}
v^-(x,y) + v^+(x,y)  \leq    v(x) + v(y). 
\end{equation} 
\end{definition}

In the sequel we shall assume, as  part of the definition of $v: P \rightarrow \mathbf{R}$ being an upper or lower  valuation,   the associated  condition on  filters or ideals in $P$.  Definition \ref{def:1}  generalizes the definitions given by Monjardet \cite{MoB81}, Leclerc \cite{LeB93}.  A benefit of considering both isotone and antitone  valuations  may be seen in Proposition \ref{prop:logs.val.new}.
\begin{definition}[Monjardet \cite{MoB81}] 
 Let $P$ be a poset with $0$.  An isotone  function 
 $v: P \rightarrow \mathbf R$ is a 
 \emph{lower valuation} if for all $x,y,z \in P$ with 
 $x \leq z$, $y \leq z$, 
 \begin{equation}
  v(x) + v(y)  \leq v^-(x,y) + v(z).
 \end{equation} 
 Let $P$ be a poset with $1$.   An isotone  function 
 $v: P \rightarrow \mathbf R$  is an 
 \emph{upper valuation} if for all $x,y,z \in P$ with 
 $ z \leq x$, $z \leq y$, 
 \begin{equation}
  v^+(x,y) + v(z)  \leq    v(x) + v(y). 
 \end{equation} 
\end{definition}
\begin{definition}[Leclerc \cite{LeB93}] Let $ L = (L, \leq, \wedge) $ be a $\wedge$-semilattice.  A strictly isotone function $v: L \rightarrow \mathbf R$ is a \emph{lower valuation}  if and only if it satisfies the following property whenever $x \vee y$ exists:
\begin{equation} 
v(x) + v(y)  \leq v(x \vee y) + v(x \wedge y) .
\end{equation} 
 Let $L = (L, \leq, \vee)$ be a $\vee$-semilattice.  
A  strictly isotone function $v: L \rightarrow \mathbf R$ is an  \emph{upper valuation} if and only if it satisfies the following property whenever $x \wedge y$ exists: 
\begin{equation} 
 v(x \vee y) + v(x \wedge y) \leq v(x) + v(y). 
\end{equation}  
\end{definition}

 \begin{proposition}
 Let   $P$ be a finite  poset   equipped with a  valuation $v(x): P \rightarrow \mathbf{R}$ having the properties listed (by row)  in Table \ref{tab:metrics}.      Then the corresponding  formula for $d_v(x,y)$ defines a metric on $P$.  
 \end{proposition}
\begin{table}[t]
 \begin{center}
  \begin{tabular}{|r|  r  c |}  \hline  
  $\forall x,y \in P,\,\,$  & 
  valuation $v(x)$ \raisebox{1ex}{\phantom{!}}  &
  metric \\[.5ex] \hline 
  \raisebox{1ex}{\phantom{!}} 
  $\down x \cap \down y \neq \varnothing$ \, &  
  \,strictly isotone, lower: \,  &  
  $d_v(x,y) = v(x) + v(y) - 2 v^-(x,y)$ \,  \\[1ex] 
  $\up x \cap \up y \neq \varnothing$ \,   &   
  \,strictly antitone, lower:  \,  & 
  $d_v(x,y) = v(x) + v(y) - 2 v^+(x,y)$ \,  \\[1ex] 
  $\up x \cap \up y \neq \varnothing$ \,  &
  \,strictly isotone, upper: \,   & 
  $d_v(x,y) = 2v^+(x,y) - v(x) - v(y)$ \,  \\[1ex] 
  $\down x \cap \down y \neq \varnothing$ \, &
   \,strictly antitone, upper: \, & 
   $d_v(x,y) =2v^-(x,y) -v(x) - v(y)$ \,  \\[.5ex]  \hline
   \end{tabular} 
   \end{center} 
\vspace{1ex} 
\caption{\label{tab:metrics} We assume $P$ is finite; then  $d_v(x,y) = 0 \, \Rightarrow \,  x= y$.}
\end{table}
{\it Proof.} Suppose that  $v$ is a strictly isotone lower valuation; the other cases  are similar.  We  verify the triangle inequality.     Fix  $x,y,z \in P$. The inequality  
 $v^-(x,y) + v^-(y,z) \leq v^-(x,z) + v(y)$,
 which we now establish,    implies  $  d_v(x,z)  \leq  d_v(x, y) + d_v(y,z)$.     Let
\begin{align}
\alpha &\in \{ p \in  \down x \cap \down y  : v(p) = v^-(x,y) \}, \\[1ex]
\beta &\in \{ p \in \down y \cap \down z   : v(p) = v^-(y,z) \}. 
\end{align}
(These sets are nonempty by the hypothesis that $P$ is finite.)  
Since $v$ is an isotone lower valuation, we have
$$v^-(x,y) + v^-(y,z) = v(\alpha) + v(\beta)  \leq v^-(\alpha, \beta) + v(y),$$
and since $\alpha  < x$, $\beta< z$, it follows that $v^-(\alpha, \beta) < v^-(x,z)$.  \qed

If the valuation is merely isotone (or antitone), then the corresponding $d_v(x,y)$ is a quasimetric, which is defined by relaxing the metric condition `$d(x,y) = 0 \, \Rightarrow \, x = y$'.   Note that if $P$ is not finite, then  $d_v(x,y)$  is  a quasimetric but it need not be a metric.   
\subsection{Bounds on $d_v(x,y)$}
Before turning to examples we note the following bounds, and a condition for their universal attainment.   
\begin{proposition} \label{prop:d_v}  
Suppose $v: P\rightarrow \mathbf{R}$  is either an upper or lower valuation. Then 
\begin{equation}\label{eq:dist.v+v-}
 \begin{cases}
  d_v(x,y) \leq v^+(x,  y) - v^-(x, y) & \text{ if $v$ is isotone,}
  \\[.5ex]
  d_v(x,y) \leq v^-(x,  y) - v^+(x, y) & \text{ if $v$ is antitone}.
 \end{cases}
\end{equation} 
In either case, equality holds for all $x,y \in P$ if and only if $v$ is both an upper and lower valuation.   
\end{proposition}
{\it  Proof.}  Both of these assertions follow directly from the definitions. 
\qed

It turns out that if $P$ is a $\vee$-semilattice with $0$, then this upper bound for $d_v(x,y)$
 is rarely attained simultaneously for all $x,y \in P$. The following  is  essentially  Proposition \ref{prop:d_v}  combined with  \cite[Theorem 3]{Ramana.Murty:Engelbert:85}.  

\begin{proposition} Let $L = (L, \leq, \vee) $ be a finite $\vee$-semilattice with $0 \in L$, and $v : L \rightarrow \mathbf{R}$ be a strictly isotone  upper or lower valuation.     Equality can not hold in (\ref{eq:dist.v+v-}) for all $x,y \in L$ unless $L$ is a modular lattice.
\end{proposition}
{\it Proof.}   Let $z_0  \in \{ z \in \down x \cap \down y  : v(z) = v^-(x,y) \} \neq \varnothing$.   
Then $z_0$ is a lower bound of $x$ and $y$.   
Let $c$ be any other lower bound of $x$ and $y$.  Then $c \vee z_0$ is a lower bound of both $x$ and $y$ which implies $v(c \vee z_0) \leq v(z_0)$. But since $v$ is isotone $v(z_0) \leq v(c \vee z_0)$.  Strict isotonicity of $v$ implies $z_0 = c \vee z_0 $, hence $z_0 \geq c$ and $z_0$ is the greatest lower bound of $x$ and $y$.  Therefore $L$ is a lattice. Now assume that equality holds in (\ref{eq:dist.v+v-}) for all $x,y \in L$.  By Proposition \ref{prop:d_v}, $v$ is both an upper and lower valuation, hence $v$ is a valuation on $L$ (meaning that $v(x) + v(y) = v(x \vee y) + v(x \wedge y)$ for all $x, y \in L$). It is well known that the existence of a strictly isotone valuation on a lattice $L$ implies that  $L$ is a modular lattice \cite{Birkhoff:79}.  
\qed  

\subsection{Examples and discussion}
On a finite $\wedge$-semilattice $v_*(x) = |\down x|$  is an isotone
lower valuation; on a finite  $\vee$-semilattice $v^*(x) = |\up
x|$ is an antitone lower valuation.    More generally we have the following, in which cardinality is replaced by a sum over a nonnegative  weighting function:  
\begin{proposition}\label{prop:a}
Let  $P$ be a finite $\wedge$-semilattice, $t(x):  P \rightarrow [0,\infty)$ a non-negative weighting function.  Then the map $v_*(x) : P \rightarrow [0,\infty)$ defined by  
\begin{equation} \label{eq:def.v.lower.star}
v_*(x) = \sum_{x' \leq x} t(x')
\end{equation} is an isotone lower valuation. If $t(x) $ is strictly positive then $v_*(x)$ is strictly isotone. 
\end{proposition}
{\it Proof.}
Since $P$ is a $\wedge$-semilattice,  $v_*^-(x,y) = v_* (x \wedge y)$; it is sufficient to establish that for all $x,y,z \in P $ such that $x \leq z,  y \leq z$, that 
  $v_*(x) + v_*(y) \leq v_*(z) + v_*(x \wedge y )$.   
Fix $x, y \in P$ and let $J_x$ and $J_y$ denote the disjoint sets 
 $J_x  =  \down x  \cap  (\down {x\wedge y}) ^c$,  $J_y  
 =   \down y  \cap  (\down {x \wedge y} )^c$.  
For any $z \in P$ such that $x \leq z$, $y \leq z$, we have the disjoint union  and inclusion:  
\begin{equation}
J_x \cup  J_y \cup   \down{x\wedge y}  = \down x \cup \down y \subseteq \down z. 
\end{equation}
Then 
\begin{align}
v_*(x) & + v_*(y) - v_*(x \wedge y) =  \sum_{w \leq x} t(w) +  \sum_{w \leq y} t(w)  - \sum_{w \leq   x \wedge y  } t(w) \\[1ex]
& = \sum_{w \in J_x} t(w) + \sum_{w \in J_y} t(w) +  \sum_{w \in \down{ x \wedge y}  } t(w)  \leq v(z). \,\qed
\end{align}
The proof of the following proposition is similar and is omitted: 
\begin{proposition}
Let  $P$ be a finite $\vee$-semilattice, $t(x)\!:\!  P \rightarrow [0,\infty)$ a non-negative weighting function.  Then the map $v^*(x) : P \rightarrow [0,\infty)$ defined by  
\begin{equation} \label{eq:def.v.upper.star}
 v^*(x) = \sum_{x \leq x'} t(x')
\end{equation} 
is an antitone lower valuation. If $t(x) $ is strictly positive then $v_*(x)$ is strictly antitone. 
\end{proposition}

If $t(x)$ is the  indicator function of any subset $K \subseteq P$ , where $P$ is a finite $\vee$-semilattice, then $v^*(x)$ as given by (\ref{eq:def.v.upper.star}), is a lower valuation and $\kappa(x) = A - v^*(x)$ is an upper valuation for any $A \in \mathbf{R}$.  Letting $K$ denote the meet-irreducible elements of $P$, $K(x) = \{ k \in K : x \leq k\}$, and $A = |K|$, yields  the  upper valuation $\kappa(x) = |K / K(x)|$   given in \cite{LeB93}. (The use of meet-irreducible elements  is not necessary for defining the upper valuation $\kappa(x)$ given in  \cite{LeB93}:  we may replace $K$ by any subset of $P$ and obtain an upper valuation.) 

If $P$ is a poset that is not a $\wedge$-semilattice, then $v_*(x) = |\down x |$   need not be a lower valuation.  For example, $v_*(x) = |\down x|$ is not a lower valuation on the poset defined by the covering relations: $0 \prec \{ a, b, c\} \prec \{ d, e\} \prec 1$. Similarly if $P$ is a poset that is not a $\vee$-semilattice then  $v^*(x) = |\up x |$ need not be a lower valuation.  

To extend this counterexample, we consider sufficient  conditions for $v_*(x) = |\down x |$ and  $v^*(x) = |\up x |$   to be lower valuations:  let ${\mathscr P^\star \hspace{-.15ex}}$ denote the collection of finite bounded partially ordered sets, which includes all finite lattices.    A measure of the degree to which a poset $P \in {\mathscr P^\star \hspace{-.15ex}}$ deviates  from being a 
$\wedge$-semilattice or $\vee$-semilattice  (which are equivalent for  $P \in {\mathscr P^\star \hspace{-.25ex}}$) is given by the functions  
$\Delta_\wedge, \Delta_\vee: \mathscr P^\star \rightarrow \mathbb N_0$, defined by  
\begin{eqnarray}
 & \quad \Delta_\wedge(P)  = \max_{x,y \in P } D_\wedge(x,y),  &D_\wedge(x,y)  = 
 | \down x \cap \down y |- 
 \max \big\{| \down z| :  z \in \mathbf{max}( \down x \cap \down y)  \big\}, \\[.8ex]
 &\quad \Delta_\vee(P) = \max_{x,y \in P } D_\vee(x,y),    & D_\vee(x,y)   = 
 | \up x \cap \up y |- 
 \min \big\{ | \up z|  : z \in \mathbf{min}(\up x \cap \up  y) \big\}.
\end{eqnarray}

\begin{proposition}
Suppose  $P \in \mathscr P^\star$.  Then $P$  is a lattice if and only if  $\Delta_\wedge(P) = 0$ or $\Delta_\vee(P) = 0 $.  If $\Delta_\wedge(P) \leq 1$  then $v_*(x) = |\down x|$ is a lower valuation on $P$.  If $\Delta_\vee(P) \leq 1$  then $v^*(x) = |\up  x|$ is a lower valuation on $P$.  
\end{proposition}  
{\it Proof.}  The first assertion follows directly  from the definitions.  For the second, assume $\Delta_\wedge(P) \leq 1$. Accordingly, 
 \begin{equation}
v_*^+(x,y) =   \max \big\{  | \down z| : z \in \mathbf{max} ( \down x \cap \down y)\big\}  \geq  |\down x \cap \down y | - 1.
\end{equation}
We also have
  \begin{align} 
v_*^-(x,y) &\geq 1 + |\down x / ( \down x \cap \down y) | 
+ | \down y / ( \down x \cap \down y) | +  |\down x \cap \down y | , \\[.5ex]
v_*(x) &= |\down x / ( \down x \cap \down y) | +  |\down x \cap \down y |, \\[.5ex]
v_*(y) &=  | \down y / ( \down x \cap \down y) | +  |\down x \cap \down y |; 
\end{align}so that  $v_*(x)$ satisfies (\ref{eq:def.lower}), and is therefore a lower valuation.   The case  $\Delta_\vee(P) \leq 1$ is similar.  \qed

Our next example leads  back to the question on logarithms:   
\begin{example} \label{ex:group.theory}
Let $G = (G, \cdot, e)$ be a multiplicative group  and  
$\mathscr L = (\mathscr L, \subseteq)$ be the collection of finite subgroups of $G$, partially ordered by inclusion.  Then 
$\mathscr L$ is a $\wedge$-semilattice in which $X \wedge Y = X  \cap Y$.   The maps  $c(X) = |X|$ and $v(X) = \log |X|$ are both lower valuations on $\mathscr L$, the latter inducing  the so-called \emph{finite subgroup metric:}
\begin{equation}
d_v(X,Y) = \log \frac{|X | |Y|}{( | X \cap Y| )^2}.
\end{equation}
  If $G$ is abelian, then $\mathscr L$ is a lattice (but $\mathscr L $ is not necessarily a complete lattice) and $v(X) = \log |X|$ is an upper valuation as well.      
\end{example}
{\it Proof.} Whether or not $XY $ is a subgroup of $G$, the product formula \cite[p.~14]{Rotman:84}  states that   
\begin{equation}\label{eq:product.formula}
|X| |Y| = |XY| |X \cap Y|. 
\end{equation}   Let $m  =[ X\!:\! X \cap Y] $ and $n =  [Y\!:\! X\cap Y]$. Then (\ref{eq:product.formula}) implies    $|XY| = mn|X\cap Y|$, and since $m + n \leq mn + 1$ for all  $m,n \in \mathbb{Z}^+$ it follows that $|X| + |Y| \leq |X \cap Y| + |X Y |$. Hence  if  $X \vee Y \in \mathscr L$, then $|X| + |Y| \leq |X \wedge Y| + |X \vee Y | $ and $c(X)$ is a lower valuation.   The fact that $v(X) = \log |X|$ is a lower valuation follows from  (\ref{eq:product.formula}) and $XY \subseteq X\vee Y$. If $G$ is abelian then $|X \vee Y| = |XY|$ and $v(X)$ is also an upper valuation.  If $G$ is an infinite abelian group then $X, Y \in \mathscr L \Rightarrow X \wedge Y \in \mathscr L$, $X\vee Y \in \mathscr L$ ($X \vee Y$ is finite) so that $\mathscr L$ is a lattice, but the the join over an arbitrary number of finite subgroups need not be finite so $\mathscr L$ need not be a complete lattice.    
\qed 
 \section{Composition with logarithms}
Suppose  $v(x):P \rightarrow \mathbf R$ is either an upper valuation or  lower  valuation, either isotone or antitone.    Observe that 
\begin{equation}v^\prime(x) =  K \cdot v(x) + A, \quad K \in \mathbf R/\{0\}, \, A \in \mathbf R,
\end{equation}
is also an upper  or lower valuation,  and if $K <0$,  \emph{upper} and  \emph{lower} are interchanged,  as well as \emph{isotone} and \emph{antitone}. 

\begin{proposition}\label{prop:logs.val}   
Suppose $u :  P \rightarrow \mathbf R^+$ is a strictly positive isotone (antitone) upper valuation.  Then $\ell(x) = \log u(x)$ is an isotone (antitone) upper valuation.  On the other hand, if $v(x) :   P \rightarrow \mathbf R^+$ is a strictly positive isotone lower valuation then $\ell^\prime(x) = \log v(x) $ need not be an upper valuation or a lower valuation. 
\end{proposition}  
{\it Proof.} 
  Let $x,y \in P$, and   let  $a = u^+(x,y)$, $b = u(x)$, $c= u(y)$, $d= u^-(x,y)$. We  treat the case that $u(x)$ is isotone. By hypothesis 
\begin{align}  \label{eq:ad.bc}
& a + d \leq b + c, \quad  a ,b, c, d > 0, \\ \label{eq:d.min.max.a}
&  d \leq \min\{ b,c\} \leq \max\{b,c\} \leq a. 
\end{align} 
 Since   $\ell(x) = \log u(x)$ is isotone, $\ell^+(x,y) = \log u^+(x,y)$ and $\ell^-(x,y) = \log u^-(x,y)$.    The function $\ell(x)$ is an upper valuation because it  satisfies (\ref{eq:def.upper}), that is,  
\begin{equation} 
\ell^+ (x,y) + \ell^-(x,y)  =   
\log a + \log d \leq \log b + \log c  = \ell(x) + \ell(y),
\end{equation} 
or equivalently, $ad \leq bc$.  Indeed, let $d = \min\{b,c\} - X$, $a = \max\{ b,c \} + Y$, where $X, Y \geq 0$.  Note that $Y \leq X$ follows from (\ref{eq:ad.bc}).  Since $bc = \min\{b, c\} \cdot \max\{ b,c \}$, we have 
\begin{align}\nonumber
ad &=  bc  + Y \cdot \min\{b,c\} - X \cdot \max\{ b,c \} - XY   \\[.1ex] 
& \leq   bc  + X   \big(\! \min\{b,c\} -  \max\{ b,c \} \big)  - XY \,  \leq  bc.
\end{align}
The case that $u(x)$ is antitone may be treated similarly.  If $u(x)$ is antitone, then instead of (\ref{eq:d.min.max.a}) we have $a \leq \min\{ b,c\} \leq \max\{b,c\} \leq d$, while (\ref{eq:ad.bc}) still holds. 

Finally, as a counterexample, consider  the lower valuation $v(x) = \sum_{x' \leq x} t(x') $ defined on the Boolean lattice $M_2$ with covering relations $0 \prec \{p, q\} \prec 1$, where  $t: M_2 \rightarrow \mathbf R^+$ is a discrete probability distribution ($v(1) = 1$).   Then $\log v(x)$ need not be an upper or  lower valuation (depending on $t(x)$). \qed

Combining Proposition \ref{prop:logs.val} with the observation  preceding its statement yields the following:  
\begin{proposition}\label{prop:logs.val.new}   
Suppose $u :\! P \rightarrow \mathbf R$ is an 
isotone (antitone) upper valuation.  Then $L(x) =\log(K \cdot u(x) + A)$ is an  isotone (antitone) upper valuation for any $K > 0$ and  $A > -\min_{x \in P} K \cdot u(x)$.  Suppose $v :\! P \rightarrow \mathbf R $ is an isotone (antitone) lower valuation.  Then $L^\prime(x) = -\log (K \cdot v(x) + A)$ is an isotone (antitone) lower  valuation for any $K < 0$ and $A >  \max_{x \in P}  | K | \cdot v(x)$.  
\end{proposition}
While the valuations $L(x)$ and $L'(x)$ of Proposition \ref{prop:logs.val.new} are available for defining metrics  on $P$,   we note  that the formula of  Jiang and Conrath \cite{JiJCoD97} 
\begin{equation}
 \text{dist}_{JC} (x, y) = 
 \mathcal  I(x) + \mathcal  I(y) - 2  \mathcal  I^+(x,y), 
 \quad  
 \mathcal  I(x) = -\log p(x),  
\end{equation}
in which $p(x)$ is a cumulative probability of the form (\ref{eq:def.v.lower.star}),  is not necessarily a metric defined on a general poset.    
As a counterexample, consider the poset defined by the covering relations:  $\{ z_1, z_2 \} \prec a$, $\{ z_1, z_3 \} \prec b$, $\{ z_2, z_3 \} \prec c$, $\{a, b,c \} \prec 1$ with discrete probability distribution $t(x)$ and cumulative probability $p(x) = \sum_{x' \leq x} t(x)$.  Then $\text{dist}_{ JC} (x, y)$ need not be a metric: depending on $t(x)$, it can happen that  $$\text{dist}_{ JC} (z_1, z_2)  + \text{dist}_{JC} (z_2, z_3) \leq  \text{dist}_{ JC} (z_1, z_3).$$  A sufficient condition for $\text{dist}_{JC} (x, y)$ to be a metric is that the poset be a tree.    
 
\bibliographystyle{plain}
 \bibliography{Poset.val}

\end{document}